\documentclass[11pt]{article}
\usepackage{amsfonts}
\usepackage{amssymb}
\usepackage{amsmath}
\usepackage{latexsym}
\usepackage{amsthm}
\usepackage[english]{babel}
\usepackage{epsfig}
\newtheorem{theorem}{Theorem}

\title{A solution to the tennis ball problem}
\author{Anna de Mier \qquad Marc Noy  \\
        \normalsize\smallskip Universitat Polit\`{e}cnica de
        Catalunya\thanks{Departament de Matem\`{a}tica Aplicada II.
        Pau Gargallo 5, 08028 Barcelona, Spain. E-mail addresses:
         {\tt \{anna.de.mier,marc.noy\}@upc.es}. Research partially supported
         by projects BFM2001-2340 and CUR Gen. Cat. 1999SGR00356.}
         }
\date{}

\begin{document}
\maketitle

\begin{abstract}
We present a complete solution to the so-called tennis ball
problem, which is equivalent to counting lattice
paths in the plane that use North and East steps and lie between
certain boundaries. The solution takes the form of explicit
expressions for the corresponding generating functions.

Our method is based on the properties of Tutte polynomials of
matroids associated to lattice paths. We also show how the same
method provides a solution to a wide generalization of the
problem.

\end{abstract}

\section{Introduction}

The statement of the tennis ball problem is the following. There
are $2n$ balls numbered $1,2,3,\dots,2n$. In the first turn balls
1 and 2 are put into a basket and one of them is removed. In the
second turn balls 3 and 4 are put into the basket and one of the
three remaining balls is removed. Next balls 5 and 6 go in and one
of the four remaining balls is removed. The game is played $n$
turns and at the end there are exactly $n$ balls outside the
basket. The question is how many different sets of balls may we
have at the end outside the basket.

It is easy to reformulate the problem in terms of lattice paths in
the plane that use steps $E=(1,0)$ and $N=(0,1)$. It amounts to
counting the number of lattice paths from $(0,0)$ to $(n,n)$ that
never go above the path $NE\cdots NE = (NE)^n$. Indeed, if $\pi =
\pi_1 \pi_2 \dots \pi_{2n-1}\pi_{2n}$ is such a path, a moment's
thought shows that we can identify the indices $i$ such that
$\pi_{2n-i+1}$ is a $N$ step with the labels of balls that end up
outside the basket. The number of such paths is well-known to be a
Catalan number, and this is the answer obtained in \cite{GM}.

The problem can be generalized as follows \cite{tennis}. We are
given positive integers $t<s$ and $sn$ labelled balls. In the
first turn balls $1,\dots,s$ go into the basket and $t$ of them
are removed. In the second turn balls $s+1,\dots,2s$ go into the
basket and $t$ among the remaining ones are removed. After $n$
turns, $tn$ balls lie outside the basket, and again the question
is how many different sets of balls may we have at the end.
Letting $k=t,l=s-t$, the problem is seen as before \underline{}to
be equivalent to counting the number of lattice paths from $(0,0)$
to $(ln,kn)$ that use $N$ and $E$ steps and never go above the
path $N^kE^l\cdots N^kE^l = (N^kE^l)^n$. This is the version of
the problem we solve in this paper.

\begin{figure}[htb]
\begin{center}
\epsfysize=5cm\epsffile{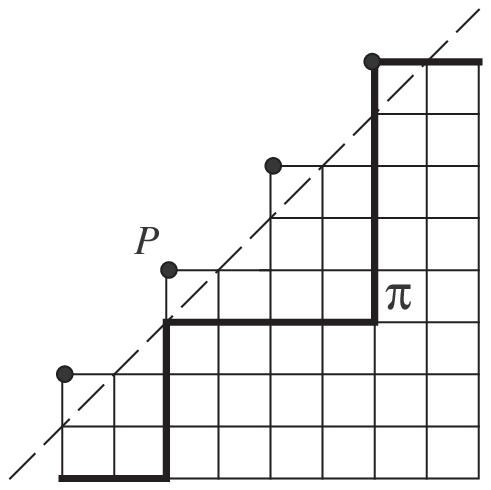} \caption{\label{fig:path}  The
path $\pi=\underline{EE}NN\underline{N}EEEENNN\underline{NN}EE$
not above $P=(N^2E^2)^4$. It has $i(\pi)=3$ and $e(\pi)=2$,
corresponding to the steps underlined.}
\end{center}
\end{figure}

From now on we concentrate on lattice paths that use $N$ and $E$
steps. To our knowledge, the only cases solved so far are $k=1$
and $k=l=2$. The case $k=1$ is straightforward, the answer being a
generalized Catalan number ${1 \over (l-1)n+1}{ln \choose n}$. The
case $k=l=2$ (corresponding to the original problem when $s=4,
t=2$) is solved in \cite{tennis} using recurrence equations; here
we include a direct solution. This case is illustrated in Fig.
\ref{fig:path}, to which we refer next. A path $\pi$ not above
$(N^2E^2)^n$ is ``almost'' a Catalan path, in the sense that it
can raise above the dashed diagonal line only through the dotted
points. But clearly between two consecutive dotted points hit by
$\pi$ we must have an $E$ step, followed by a Catalan path of odd
semilength, followed by a $N$ step. Thus, $\pi$ is essentially a
sequence of Catalan paths of odd semilength. If $G(z) = \sum_n {1
\over n+1}{2n \choose n}z^n$ is the generating function for the
Catalan numbers, take the odd part $G_o(z) = (G(z)-G(-z))/2$. Then
expand $1/(1-zG_o(z))$ to obtain the sequence
$1,6,53,554,6363,\dots$, which agrees with the results in
\cite{tennis}.

Let $P$ be a lattice path from $(0,0)$ to $(m,r)$, and let $b(P)$
be the number of paths from $(0,0)$ to $(m,r)$ that never go above
$P$. If $P N$ denotes the path obtained from $P$ by adding a $N$
step at the end of $P$, then clearly $b(P) = b(PN)$. However, it
is not possible to express $b(PE)$ simply in terms of $b(P)$,
where $PE$ has the obvious meaning. As is often the case in
counting problems, one has to enrich the objects under enumeration
with additional parameters that allow suitable recursive
decompositions. This is precisely what is done here: equations
(\ref{eq:t1}) and (\ref{eq:t2}) in the next section contain
variables $x$ and $y$, corresponding to two
parameters that we define on lattice paths not above a given path
$P$. These equations are the key to our solution.

The basis of our approach is the connection between lattice paths
and matroids established in \cite{BMN}, where the link with the
tennis ball problem was already remarked. For completeness, we
recall the basic facts needed from \cite{BMN} in the next section.
In Section \ref{sec:main} we present our solution to the tennis
ball problem, in the form of explicit expressions for the
corresponding generating functions; see Theorem~\ref{th:main}. In
Section \ref{sec:gen} we show how the same method can be applied
to a more general problem. We conclude with some remarks.

\section{Preliminaries}

The contents of this section are taken mainly from \cite{BMN},
where the reader can find additional background and references on
matroids, Tutte polynomials and lattice path enumeration.

A \emph{matroid} is a pair $(E, \mathcal{B})$ consisting of a
finite set $E$ and a nonempty collection $\mathcal{B}$ of subsets
of $E$, called \emph{bases} of the matroid, that satisfy the
following conditions: (1) No set in $\mathcal{B}$ properly
contains another set in $\mathcal{B}$, and (2) for each pair of
distinct sets $B,B'$ in $\mathcal{B}$ and for each element $x\in
B-B'$, there is an element $y\in B'-B$ such that $(B-x)\cup y$ is
in $\mathcal{B}$.

Let $P$ be a lattice path from $(0,0)$ to $(m,r)$. Associated to
$P$ there is a matroid $M[P]$ on the set $\{1,2,\dots,m+r\}$ whose
bases are in one-to-one correspondence with the paths from $(0,0)$
to $(m,r)$ that never go above $P$. Given such a path
$\pi=\pi_1\pi_2\dots\pi_{m+r}$, the basis corresponding to $\pi$
consists of the indices $i$ such that $\pi_i$ is a $N$ step.
Hence, counting bases of $M[P]$ is the same as counting lattice
paths that never go above $P$.

For any matroid $M$ there is a two-variable polynomial with
non-negative integer coefficients, the Tutte polynomial
$t(M;x,y)$. It was introduced by Tutte \cite{tutte2} and presently
plays an important role in combinatorics and related areas (see
\cite{welsh}). The key property in this context is that $t(M;1,1)$
equals the number of bases of $M$.

Given a path $P$ as above, there is a direct combinatorial
interpretation of the coefficients of $t(M[P];x,y)$. For a path
$\pi$ not above $P$, let $i(\pi)$ be the number of $N$ steps that
$\pi$ has in common with $P$, and let $e(\pi)$ be the number of
$E$ steps of $\pi$ before the first $N$ step, which is 0 if $\pi$
starts with a $N$ step. See Fig. \ref{fig:path} for an
illustration.

Then we have (see \cite[Th. 5.4]{BMN})
\begin{equation}\label{eq:act}
    t(M[P];x,y) = \sum_{\pi} x^{i(\pi)} y^{e(\pi)},
\end{equation}
where the sum is over all paths $\pi$ not above $P$. A direct
consequence is that $t(M[P];1,1)$ is the number of such paths.

Furthermore, for the matroids $M[P]$ there is a rule for computing
the Tutte polynomial that we use repeatedly (see \cite[Section
6]{BMN}). If $P N$ and $P E$ denote the paths obtained from $P$ by
adding a $N$ step and an $E$ step at the end of $P$, respectively,
then
\begin{eqnarray}
\label{eq:t1} t(M[P N];x,y) &=& x\,t(M[P],x,y), \\
\label{eq:t2} t(M[P E];x,y) &=& {x\over x-1} \, t(M[P],x,y) +
\left(y-{x \over x-1} \right) t(M[P];1,y).
\end{eqnarray}
The right-hand side of (\ref{eq:t2}) is actually a polynomial,
since $x-1$ divides $t(M[P];x,y)-t(M[P];1,y)$. The key observation
here is that we cannot simply set $x=y=1$ in (\ref{eq:t2}) to
obtain an equation linking $t(M[PE];1,1)$ and $t(M[P];1,1)$.

For those familiar with matroid theory, we remark that $M[PN]$ and
$M[PE]$ are obtained from $M[P]$, respectively, by adding an
isthmus and taking a free extension; it is known that formulas
(\ref{eq:t1}) and (\ref{eq:t2}) correspond precisely to the effect
these two operations have on the Tutte polynomial of an arbitrary
matroid.

From (\ref{eq:act}) and the definition of $i(\pi)$ and $e(\pi)$,
equation (\ref{eq:t1}) is clear, since any path associated to
$M[PN]$ has to use the last $N$ step. For completeness, we include
a direct proof of equation (\ref{eq:t2}).

We first rewrite the right-hand side of (\ref{eq:t2}) as

\begin{eqnarray*}
& &\frac{x}{x-1}(t(M[P];x,y)-t(M[P];1,y)) + yt(M[P];1,y)=\\
& &\sum_{\pi}\frac{x}{x-1}y^{e(\pi)}(x^{i(\pi)}-1)+y^{e(\pi)+1} =\\
& &\sum_{\pi} y^{e(\pi)}(y+x+x^2+\cdots+x^{i(\pi)}),\\
\end{eqnarray*}
where the sums are taken over all paths $\pi$ that do not go above
$P$.

\begin{figure}[htb]
\begin{center}
\epsfysize=7cm\epsffile{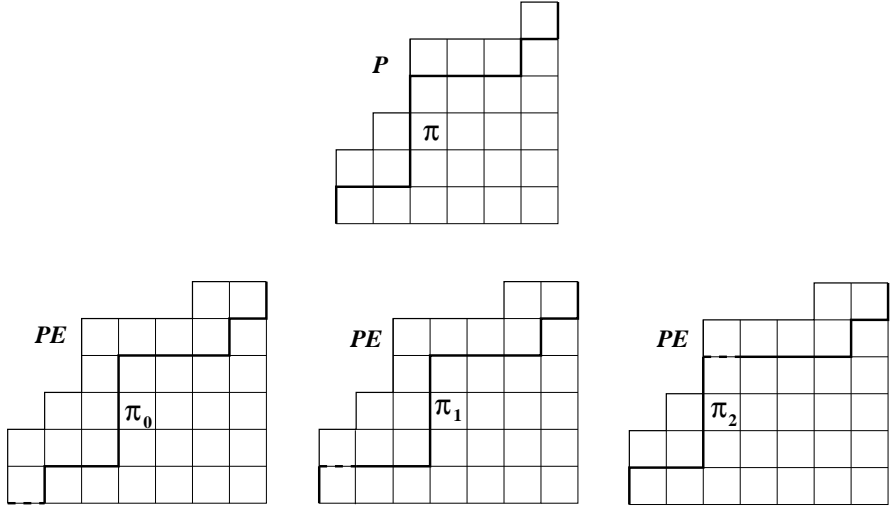}
\caption{\label{fig:activities}Illustrating the combinatorial
proof of formula (\ref{eq:t2}).}
\end{center}
\end{figure}

To prove the formula, for each path $\pi$ not above $P$ we find
$i(\pi)+1$ paths not above $PE$ such that their contribution to
$t(M[PE];x,y)$ is $y^{e(\pi)}(y+x+x^2+\cdots+x^{i(\pi)})$.
Consider first the path $\pi_0=E\pi$; it clearly does not go above
$PE$ and its contribution to the Tutte polynomial is
$y^{e(\pi)+1}$. Now for each $j$ with $1\leq j \leq i(\pi)$,
define the path $\pi_j$ as the path obtained from $\pi$ by
inserting an $E$ step after the $j$th $N$ step that $\pi$ has in
common with $P$ (see Fig.~\ref{fig:activities}). The path $\pi_j$
has exactly $j$ $N$ steps in common with $PE$, and begins with
$e(\pi)$ $E$ steps.

It remains only to show that each contribution to the Tutte
polynomial of $M[PE]$ arises as described above. Let $\pi'$ be a
path that never goes above $PE$ and consider the last $N$ step
that $\pi'$ has in common with $PE$; clearly the next step must be
$E$. Let $\widetilde{\pi}$ be the path obtained after removing
this $E$ step. Since $\pi'$ had no $N$ steps in common with $PE$
after the $E$ step removed, the path $\widetilde{\pi}$ does not go
above $P$. Thus the path $\pi'$ can be obtained from
$\widetilde{\pi}$ by adding an $E$ step after the $i(\pi')$-th $N$
step that $\widetilde{\pi}$ has in common with $P$, and hence
$\pi'$ arises from $\widetilde{\pi}$ as in the previous paragraph.
It is easy to show that $\pi'$ cannot be obtained in any other way
by applying the procedure described above, and this finishes the
proof.

\section{Main result}\label{sec:main}

Let $k,l$ be fixed positive integers, and let $P_n = (N^k E^l)^n$.
Our goal is to count the number of lattice paths from $(0,0)$ to
$(ln,kn)$ that never go above $P_n$. From the considerations in
the previous section, this is the same as computing
$t(M[P_n];1,1)$. Let
 $$
    A_n = A_n(x,y) = t(M[P_n];x,y).
 $$
By convention, $P_0$ is the empty path and $A_0=1$.

In order to simplify the notation we introduce the following
operator $\Phi$ on two-variable polynomials:
 $$
    \Phi A(x,y) = {x \over x-1} A(x,y) + \left(y-{x
    \over x-1} \right) A(1,y).
 $$
Then, by equations (\ref{eq:t1}) and (\ref{eq:t2}) we have
 $$
      A_{n+1} = \Phi^l (x^k A_n),
 $$
where $\Phi^i$ denotes the operator $\Phi$ applied $i$ times.

For each $n\ge0$ and $i=1,\dots,l$, we define polynomials
$B_{i,n}(x,y)$ and $C_{i,n}(y)$ as
\begin{eqnarray*}
    B_{i,n} &=& \Phi^i\left(x^k A_n(x,y)\right), \\
    C_{i,n} &=& B_{i,n}(1,y).
\end{eqnarray*}
We also set $C_{0,n}(y) = A_n(1,y)$. Notice that
$B_{l,n}=A_{n+1}$, and $C_{0,n}(1) = A_n(1,1)$ is the quantity we
wish to compute.

Then, by the definition of $\Phi$, we have:
 $$
\begin{array}{llll}
B_{1,n} &=& {x \over x-1}x^k A_n  &+ \left(y-{x
    \over x-1} \right) C_{0,n}; \\
B_{2,n} &=& {x \over x-1} B_{1,n}  &+ \left(y-{x
    \over x-1} \right) C_{1,n}; \\
    &\cdots& \\
B_{l,n} &=& {x \over x-1} B_{l-1,n}  &+ \left(y-{x
    \over x-1} \right) C_{l-1,n}; \\
A_{n+1} &=& B_{l,n}.
\end{array}
 $$

In order to solve these equations, we introduce the following
generating functions in the variable $z$ (but recall the
coefficients are polynomials in $x$ and $y$):
 $$
    A = \sum_{n\ge0} A_n z^n, \qquad
    C_i = \sum_{n\ge0} C_{i,n} z^n, \quad  i=0,\dots,l.
$$
We start from the last equation $A_{n+1} = B_{l,n}$ and substitute
repeatedly the value of $B_{i,n}$ from the previous equation.
Taking into account that $\sum_n A_{n+1}z^n = (A-1)/z$, a simple
computation yields
 $$
    {A -1 \over z} = {x^{k+l} \over (x-1)^l} A +
    (yx-y-x)\sum_{i=1}^l  {x^{i-1} \over (x-1)^i} \, C_{l-i}.
 $$
We now set $y=1$ and obtain
\begin{equation}\label{eq:main}
 A \left((x-1)^l-zx^{k+l}\right) = (x-1)^l -
 z\sum_{i=1}^l x^{i-1} (x-1)^{l-i} \, C_{l-i},
\end{equation}
where it is understood that from now on we have set $y=1$ in the
series $A$ and $C_i$.

By Puiseux's theorem (see \cite[Chap. 6]{Stan2}), the algebraic
equation in $w$
\begin{equation}\label{eq:alg}
    (w-1)^l-zw^{k+l} = 0
\end{equation}
has $k+l$ solutions in the field $\mathbb{C}^{\,\rm fra}((z)) =
\{\sum_{n\ge n_0} a_n z^{n/N} \} $ of fractional Laurent series.
Proposition 6.1.8 in \cite{Stan2} tells us that exactly $l$ of
them are fractional power series (without negative powers of $z$);
let them be $w_1(z),\dots,w_l(z)$.

We substitute $x=w_j$  in (\ref{eq:main}) for $j=1,\dots,l$, so
that the left-hand side vanishes, and obtain a system of $l$
linear equations in $C_0,C_1,\dots,C_{l-1}$, whose coefficients
are expressions in the $w_j$, namely
\begin{equation}\label{eq:sys}
     \sum_{i=1}^l w_j^{i-1} (w_j-1)^{l-i}  zC_{l-i} =
     (w_j-1)^l, \qquad j=1,\dots,l.
\end{equation}
Notice that, since the $C_i$ are (ordinary) power series, the
solutions $w$ of (\ref{eq:alg}) that we substitute in
(\ref{eq:main}) cannot have negative powers of $z$, hence they
must be $w_1,\dots,w_l$. We remark the similarity of this
technique with the one devised by Tutte for counting rooted planar
maps (see, for instance, \cite{tutte}).

It remains only to solve (\ref{eq:sys}) to obtain the desired
series $C_0 = \sum_n A_n(1,1)z^n$. The system (\ref{eq:sys}) can
we written as
$$
    \sum_{i=0}^{l-1} \left(\frac{w_j}{w_j-1} \right)^{i} zC_{l-i-1}=
    w_j-1, \qquad j=1,\ldots,l.
 $$
The left-hand sides of the previous equations can be viewed as the
result of evaluating the polynomial $\sum_{i=0}^{l-1} (zC_{l-i-1})
X^{i}$ of degree $l-1$ at $X=w_j/(w_j-1)$, for $j$ with $1\leq j
\leq l$. Using Lagrange's interpolation formulas, we get that the
coefficient of $X^{l-1}$ in this polynomial is

$$
    zC_0 = \sum_{j=1}^{l} \frac{w_j-1}{\prod_{i\neq j}
    \left(\frac{w_j}{w_j-1}-\frac{w_i}{w_i-1} \right)}.
$$
By straightforward manipulation this last expression is equal to
$$-\prod_{j=1}^l (1-w_j) \sum_{j=1}^l
\frac{(w_j-1)^{l-1}}{\prod_{i\neq j}(w_j-w_i)}=-\prod_{j=1}^l
(1-w_j),$$ where the last equality follows from  an identity on
symmetric functions (set $r=0$ in Exercise~7.4 in~\cite{Stan2}).

Thus we have proved the following result.

\begin{theorem}\label{th:main}
Let $k,l$ be positive integers. Let $q_n$ be the number of lattice
paths from $(0,0)$ to $(ln,kn)$ that never go above the path $(N^k
E^l)^n$, and let $w_1,\dots,w_l$ be the unique solutions of the
equation
$$(w-1)^l-zw^{k+l} = 0$$
that are fractional power series. Then the generating function
$Q(z) = \sum_{n\ge0} q_n z^n$ is given by
 $$
    Q(z) = {-1 \over z} (1-w_1)\cdots(1-w_l).
 $$
\end{theorem}

\medskip

Remark that, by symmetry, the number of paths not above
$(N^lE^k)^n$ must be the same as in Theorem \ref{th:main},
although the algebraic functions involved in the solution are
roots of a different equation.

In the particular case $k=l$ the solution can be expressed
directly in terms of the generating function $G(z) = \sum_n {1
\over n+1}{2n \choose n}z^n$ for the Catalan numbers, which
satisfies the quadratic equation $G(z)=1+zG(z)^2$. Indeed,
(\ref{eq:alg}) can be rewritten as
 $$
    w = 1 + z^{1/k} w^2,
 $$
whose (fractional) power series solutions are $G(\zeta^j
z^{1/k})$, $j=0,\dots,k-1$, where $\zeta$ is a primitive $k$-th
root of unity. For instance, for $k=l=3$ (corresponding to $s=6,
t=3$ in the original problem), $\zeta=\exp(2\pi i/3)$ and we
obtain the solution
\begin{eqnarray*}
     && {-1 \over z} (1-G(z^{1/3}))(1-G(\zeta z^{1/3}))(1-G(\zeta^2
    z^{1/3})) = \\
    && \qquad\qquad\qquad 1 + 20z + 662z^2 + 26780z^3 + 1205961z^4
    +58050204 z^5 + \cdots.
\end{eqnarray*}

In the same way, if $l$ divides $k$ and we set $p=(k+l)/l$, the
solution can be expressed in terms of the generating function
$\sum_n {1 \over (p-1)n+1}{pn \choose n}z^n$ for generalized
Catalan numbers; the details are left to the reader. As an
example, for $k=4$, $l=2$, we obtain the series
\begin{eqnarray*}
     && {-1 \over z} (1-H(z^{1/2}))(1-H(-z^{1/2})) = \\
    && \qquad\qquad\qquad 1 + 15z + 360z^2 + 10463z^3 + 337269z^4 + 11599668z^5 +
    \cdots,
\end{eqnarray*}
where $H(z) = \sum_n {1 \over 2n+1}{3n \choose n}$ satisfies
$H(z)=1+zH(z)^3$.

\section{A further generalization}\label{sec:gen}

In this section we solve a further generalization of the tennis
ball problem. Given fixed positive integers
$s_1,t_1,\dots,s_r,t_r$ with $t_i < s_i$ for all $i$, let $s =
\sum s_i, t = \sum t_i$. There are $sn$ labelled balls. In the
first turn we do the following: balls $1,\dots,s_1$ go into the
basket and $t_1$ of them are removed; then balls
$s_1+1,\dots,s_1+s_2$ go into the basket and among the remaining
ones $t_2$ are removed; this goes on until we introduce balls
$s-s_r+1,\dots,s$, and remove $t_r$ balls. After $n$ turns there
are $tn$ balls outside the basket and the question is again how
many different sets of $tn$ balls may we have at the end.

The equivalent path counting problem is: given
$k_1,l_1,\dots,k_r,l_r$ positive integers with $k = \sum k_i$,
$l=\sum l_i$, count the number of lattice paths from $(0,0)$ to
$(ln,kn)$ that never go above the path $P_n = (N^{k_1}E^{l_1}
\cdots N^{k_r}E^{l_r})^n$. The solution parallels the one
presented in Section \ref{sec:main}. We keep the notations and let
$A_n = t(M[P_n];x,y)$, so that
 $$
    A_{n+1} = \Phi^{l_r} (x^{k_r} \cdots \Phi^{l_1} (x^{k_1}
    A_n)\cdots).
 $$
As before, we introduce $l$ polynomials $B_{i,n}(x,y)$ and
$C_{i,n}(y) = B_{i,n}(1,y)$, but the definition here is a bit more
involved:

\begin{equation}\label{eq:gen}
 \begin{array}{lll}
 B_{i,n} & = \Phi^i(x^{k_1} A_n), &\quad i=1,\dots,l_1; \\
 B_{l_1+i,n} & = \Phi^i(x^{k_2} B_{l_1,n}), &\quad i=1,\dots,l_2; \\
 B_{l_1+l_2+i,n} & = \Phi^i(x^{k_3} B_{l_1+l_2,n}), &\quad i=1,\dots,l_3; \\
 & \cdots \\
 B_{l-l_r+i,n} & = \Phi^i(x^{k_r} B_{l-l_r,n}), &\quad i=1,\dots,l_r. \\
\end{array}
\end{equation}
We also set $C_{0,n}(y) = A_n(1,y)$. Again, from the definition of
$\Phi$, we obtain a set of equations involving $A_n$,
$A_{n+1}=B_{l,n}$, the $B_{i,n}$ and $C_{i,n}$. We define
generating functions $A$ and $C_i$ ($i=0,\dots,l$) as in Section
\ref{sec:main}.

Starting with $A_{n+1}=B_{l,n}$, we substitute repeatedly the
values of the $B_{i,n}$ from previous equations and set $y=1$.
After a simple computation we arrive at
\begin{equation}\label{eq:gen2}
 A \left((x-1)^l-zx^{k+l}\right) = (x-1)^l + z\,U(x,C_0,\dots,C_{l-1}),
\end{equation}
where $U$ is a polynomial in the variables $x,C_0,\dots,C_{l-1}$.
Observe that the difference between (\ref{eq:gen2}) and equation
(\ref{eq:main}) is that now $U$ is not a concrete expression but a
certain polynomial that depends on the particular values of the
$k_i$ and $l_i$.

Let $w_1,\dots,w_l$ be again the power series solutions of
(\ref{eq:alg}). Substituting $x=w_j$ in (\ref{eq:gen2}) for
$j=1,\dots,l$, we obtain a system of linear equations in the
$C_i$. Since the coefficients are rational functions in the $w_j$,
the solution consists also of rational functions; they are
necessarily symmetric since the $w_j$, being conjugate roots of
the same algebraic equation, are indistinguishable.

Thus we have proved the following result.

\begin{theorem}\label{th:maingen}
Let $k_1,l_1,\dots,k_r,l_r$ be positive integers, and let $k =
\sum k_i$, $l=\sum l_i$. Let $q_n$ be the number of lattice paths
from $(0,0)$ to $(ln,kn)$ that never go above the path
$(N^{k_1}E^{l_1} \cdots N^{k_r}E^{l_r})^n$, and let
$w_1,\dots,w_l$ be the unique solutions of the equation
$$(w-1)^l-zw^{k+l} = 0$$
that are fractional power series. Then the generating function
$Q(z) = \sum_{n\ge0} q_n z^n$ is given by
 $$
    Q(z) = {1 \over z} R(w_1,\dots,w_l),
 $$
where $R$ is a computable symmetric rational function of
$w_1,\dots,w_l$.
\end{theorem}

As an example, let $r=2$ and $(k_1,l_1,k_2,l_2)=(2,2,1,1)$, so
that $k=l=3$. Solving the corresponding linear system we obtain
 $$
    R = {(1-w_1)(1-w_2)(1-w_3) \over 2w_1 w_2 w_3 - (w_1w_2 +
    w_1w_3 + w_2w_3)},
 $$
and
 $$
 Q(z) = {1 \over z}R = 1 + 16z + 503z^2 + 19904z^3 + 885500z^4 + 42298944z^5 +
 \cdots.
 $$
It should be clear that for any values of the $k_i$ and $l_i$ the
rational function $R$ can be computed effectively. In fact, a
simple computer program could be written that on input
$k_1,l_1,\dots,k_r,l_r$, outputs $R$.

\section{Concluding Remarks}

It is possible to obtain an expression for the generating function
of the full Tutte polynomials $A_n(x,y)$ defined in
Section~\ref{sec:main}. We  have to find the values of
$C_0,C_1,\dots,C_{l-1}$ satisfying the system (\ref{eq:sys}) and
substitute back into (\ref{eq:main}). After some algebraic
manipulation, the final expression becomes
$$
    \sum_{n\ge0} A_n(x,y) z^n =
\frac{-(x-w_1)\cdots (x-w_l)}{(zx^{k+l}-(x-1)^l)(y-w_1(y-1))\cdots
(y-w_l(y-1))}.
$$
Taking $x=y=1$ we recover the formula stated in Theorem
\ref{th:main}.

On the other hand, references \cite{MS} and \cite{tennis} also
study a different question on the tennis ball problem, namely to
compute the sum of the labels of the balls outside the basket for
all possible configurations. For a given lattice path $P_n$, this
amounts to computing the sum of all elements in all bases of the
matroid $M[P_n]$. We remark that this quantity does not appear to
be computable from the corresponding Tutte polynomials alone.

Finally, as already mentioned, the technique of forcing an
expression to vanish by substituting  algebraic functions was
introduced by Tutte in his landmark papers on the enumeration of
planar maps. Thus the present paper draws in more than one way on
the work of the late William Tutte.

\end{document}